\documentclass[12pt]{amsart}

\newcommand{\eps}{\epsilon}

\newfont{\fnt}{cmr10 scaled 550}
\renewcommand{\eps}{\varepsilon} 

\newtheorem{theorem}{Theorem}[section]
\newtheorem{conj}{Conjecture}
\newtheorem{lemma}{Lemma}[section]
\newtheorem{corollary}{Corollary}[section]
\newtheorem{prop}{Proposition}[section]

\theoremstyle{remark}
\newtheorem{remark}{Remark}[section] 
\numberwithin{equation}{section}

\font\strange=msbm10

\renewcommand{\epsilon}{\varepsilon}

\renewcommand{\Sigma}{\varSigma}
\newcommand{\R}{{{\mathchoice  {\hbox{$\textstyle{\text{\strange R}}$}}
{\hbox{$\textstyle{\text{\strange R}}$}}
{\hbox{$\scriptstyle  N\kern-0.3em  R$}}  
{\hbox{$\scriptscriptstyle  R\kern-0.2em  R$}}}}}

\newcommand{\Z}{{{\mathchoice  {\hbox{$\textstyle{\text{\strange Z}}$}}
{\hbox{$\textstyle{\text{\strange Z}}$}}
{\hbox{$\scriptstyle  Z\kern-0.3em  Z$}}
{\hbox{$\scriptscriptstyle  Z\kern-0.2em  Z$}}}}}

\newcommand{\N}{{{\mathchoice  {\hbox{$\textstyle{\text{\strange N}}$}}
{\hbox{$\textstyle{\text{\strange N}}$}}
{\hbox{$\scriptstyle  N\kern-0.3em  N$}}
{\hbox{$\scriptscriptstyle  N\kern-0.2em  N$}}}}}

\renewcommand{\phi}{\varphi}

\usepackage{amsmath,amsthm,amscd}

\begin{document}
\title[The DDVV Conjecture]{Recent developments of the DDVV Conjecture}
\date{July 28, 2007}

 \author{Zhiqin Lu}
 \thanks{The 
author is partially supported by  NSF Career award DMS-0347033 and the
Alfred P. Sloan Research Fellowship.}

\address{Department of
Mathematics, University of California,
Irvine, Irvine, CA 92697}

\email[Zhiqin Lu]{zlu@math.uci.edu}

\maketitle 

\tableofcontents
\pagestyle{myheadings}

\newcommand{\nn}{{\bf n}}
\newcommand{\ka}{K\"ahler }
\newcommand{\ii}{\sqrt{-1}}

\section{Introduction}
Let $M^n$ be an immersed submanifold of $N^{n+m}(c)$, the space form of constant sectional curvature $c$. The scalar curvature of the induced metric is defined as
\[
\rho=\frac{2}{n(n-1)}\sum_{1=i<j}^n R(e_i, e_j,e_j,e_i),
\]
where $R$ is the curvature tensor of $M$ and $\{e_i\}$ is the orthonormal basis of the tangent bundle of $M$. Using the Gauss equation, in~\cite{ch1,bogd2}, it was  proved that
 \[
 \rho\leq |{\bf  H}|^2+c,
 \]
where ${\bf H}$ is the mean curvature vector. In~\cite{ddvv},  the following so-called {\it normal scalar curvature} was defined:
 \[
 \rho^{\perp}=\frac{2}{n(n-1)}\left(\sum_{1=i<j}^n\sum_{1=r<s}^m\langle
R^\perp(e_i,e_j)\xi_r,\xi_s\rangle^2\right)^{\frac 12},
\]
where $R^\perp$ is the curvature tensor of the normal bundle; $\{e_i\}$ is the orthonormal basis of the tangent bundle; and $\{\xi_j\}$ is the orthonormal basis of the normal bundle.

In the study of submanifold theory, De Smet, Dillen, Verstraelen, and Vrancken~\cite{ddvv} made the following {\sl DDVV Conjecture}: \footnote{It is also called {\it normal scalar curvature conjecture}.}

\begin{conj}\label{cc1} Let $h$ be the second fundamental form, and 
let $H=\frac 1n \,{\rm trace}\, h$ be the mean curvature tensor.  Then 
\[
\rho+\rho^\perp\leq |H|^2+c.
\]
\end{conj}

We observe that the above inequality is pointwise. Thus it is possible to rewrite the conjecture into a purely linear algebraic inequality: let $x\in M$ be a fixed point and let $(h_{ij}^r)$ ($i,j=1,\cdots,n$ and $r=1,\cdots,m$) be the coefficients of the second fundamental form under some orthonormal basis. Then by
 Suceav{\u{a}}~\cite{bogd}, or Dillen-Fastenakels-Veken~\cite{dfv}, Conjecture~\ref{cc1} can be formulated as an inequality with respect to the coefficients $h_{ij}^r$ as follows:
  \begin{align}\label{1a}
\begin{split}
&\sum_{r=1}^m\sum_{1=i<j}^n(h_{ii}^r-h_{jj}^r)^2+2n\sum_{r=1}^m\sum_{1=i<j}^n(h_{ij}^r)^2\\
&\geq 2n\left(\sum_{1=r<s}^m\sum_{1=i<j}^n\left(\sum_{k=1}^n(h_{ik}^rh_{jk}^s-h_{ik}^sh_{jk}^r)\right)^2\right)^{\frac 12}.
\end{split}
\end{align}

One can further formulate the inequality in terms of matrix notations~\cite[Theorem 2]{dfv}:

\begin{conj}[DDVV]\label{cc2}
Let $A_1,\cdots,A_m$ be symmetric $n\times n$ matrices. Then we have
\begin{equation}\label{cc2-1}
\left(\sum_{r=1}^m||A_r||^2\right)^2\geq 2\sum_{r<s}||[A_r,A_s]||^2.
\end{equation}
\end{conj}

Since the above inequality depends on $n,m$, we call it $P(n,m)$. The following special cases were known: $P(2,m)$ and $P(n,2)$ were proved in~\cite{ddvv}; $P(3,m)$ was proved in~\cite{cl}; and $P(n,3)$ was proved in~\cite{lu21}. In~\cite{dfv}, a weaker version of $P(n,m)$ was proved by using an algebraic inequality in~\cite{lianmin}. In the same paper, $P(n,m)$ was proved under the addition assumption that the submanifold is either {\it Lagrangian $H$-umbilical}, or {\it ultra-minimal} in $\mathbb C^4$.

Finally, \S 29 of the book~\cite{bychen} is a useful reference of the subject.

In this paper, we give a survey of the recent developments of the conjecture, as well as its relation to calibrated geometry, theory of random matrices, and pinching theorems of minimal submanifolds of the unit sphere. In the last part of this paper, we sketch the proofs of two important special cases: $P(3,m)$ and $P(n,3)$.

{\bf Acknowledgment.} We thank B. Suceav{\u{a}}  for  bring  the work ~\cite{ddvv} to the author's attention, for the invitation to the conference {\it Riemannian Geometry and Applications} in Bra{\c s}ov, Romania, June 2007, and for many stimulating discussions. We also thank Jason Waller for  the numerical confirmation of the conjecture.

\section{Relation to the comass problem in calibrated geometry}
Before making further analysis of Conjecture~\ref{cc2}, we recall the concept of  the comass of a $p$ form in Calibrated Geometry (cf. ~\cite{h-l}). 
 
 Consider Euclidean space $\R^n$ with orthogonal basis $e_1,\cdots,e_n$ and dual basis $dx_i=e_i^*$. Let $I=(i_1,\cdots,i_p)$ denote a multi-index with $i_1<\cdots<i_p$. Let
 \[
 \phi=\sum a_I e^*_{i_1}\wedge\cdots\wedge e_{i_p}^*
 \]
 be a $p$-covector (constant-coefficient $p$-form). 
 The comass $||\phi||^*$ of $\phi$ is given by 
 \[
 ||\phi||^*={\rm max} \{\phi(\xi)\mid\xi \text { is a $p$-plane}\}.
 \]
  For a differential form on a Riemannian manifold $M$, its comass $||\phi||^*$ is given by
 \[
 ||\phi||^*=\sup_x\,\{||\phi_x||^*\mid x\in M\}.
 \]

 In ~\cite{gmm}, Gluck, Mackenzie, and Morgan initiated the study of the comass of the
 first  Pontryagin  form on Grassmann manifolds. Later Gu~\cite{guw} generalized the results. Their results are listed as follows:
 
 \begin{theorem} The comass of the first Pontryagin form $\phi$ on the Grassmann manifold $G(n,m)$ is as follows:
 \begin{enumerate}
 \item $||\phi||^*$ is $\sqrt{3/2}$ for $n=3, m=6$, $4/3$ for $n=3, m\geq 7$, and $3/2$ for $n=4, m\geq 8$~\cite{gmm};
 \item $||\phi||^*$ is $3/2$ for $n\geq 4$ or $m\geq 8$~\cite{guw}.
 \end{enumerate}
 \end{theorem}

The definition of the comass, in the context of the comass of the first Pontryagin form, can be formulated as the following linear algebraic problem:

Let $A,B$ be two $m\times n$ matrices. Define
\[
\{AB\}=AB^T-BA^T.
\]
Let
\begin{align}
\begin{split}
&\qquad \phi(A_1\wedge A_2\wedge A_3\wedge A_4)\\
&=-\frac 12 
{\rm tr}\,(\{A_1A_2\}\{A_3A_4\}+\{A_3A_1\}\{A_2A_4\}+\{A_2A_3\}\{A_1A_4\})
\end{split}
\end{align}
for $m\times n$ matrices $A_1,A_2,A_3$, and $A_4$. The comass of $\phi$  is defined to be the maximum of the right-hand side of the  above under the condition that $A_1,A_2,A_3$, and $A_4$ are orthonormal. 

Conjecture~\ref{cc2} is similar to the above comass problem in that both problems are related to the commutator of matrices. In fact, $P(n,3)$ can be reformulated as follows: let $A,B,C$ be $n\times n$ symmetric matrices such that
\[
||A||^2+||B||^2+||C||^2=1.
\]
Then
\[
||[A,B]||^2+||[B,C]||^2+||[C,A]||^2\leq \frac 12.
\]

The major difference between these  two problems is that they have different invariant groups. The comass problem is invariant under $O(n)\times O(4m)$ (see~\cite{gmm} or ~\cite{guw} for details). On the other  hand the invariant group of the DDVV conjecture is much smaller (see \S ~\ref{p0} for details).

\section{Relation to a conjecture of B\"ottcher and Wenzel}

In~\cite{bw}, B\"ottcher and Wenzel studied the size of the commutator of two matrices $X,Y$. They showed that, if $X,Y$ are random matrices, then $||[X,Y]||^2$ is about the size $\frac 2n||X||^2\cdot||Y||^2$, which is quite small if $n$ is large. However, for fixed matrices $X,Y$, it seems that following  inequality is optimal:

\begin{conj}[B{\"o}ttcher and Wenzel]\label{cc3}
Let $X,Y$ be two $n\times n$ matrices. Then we have
\[
||[X,Y]||^2\leq 2||X||^2\cdot||Y||^2.
\]
\end{conj}

We call the above inequality $Q(n)$.
A weaker version of the  conjecture was proved~\cite[\S 3]{bw}:
\[
||[X,Y]||^2\leq 3||X||^2\cdot||Y||^2.
\]
Besides the above result, the  conjecture was proved if $X$ is of rank $1$, or if $X$ is a {\it normal} matrix in the same paper.

In this section, we give a relationship of the above BW conjecture to the DDVV conjecture. We first make the following:

\begin{conj}\label{cc4}
Let $A_1,\cdots,A_{m_1}$ be symmetric $n\times n$ matrices and let $A_{m_1+1},\cdots,A_{m_1+m_2}$ be skew-symmetric matrices. Then we have
\[
\left(\sum_{r=1}^{m_1+m_2}||A_r||^2\right)^2\geq 2\sum_{r<s}||[A_r,A_s]||^2.
\]
\end{conj}

We name the above inequality to be $P(n,m_1,m_2)$.
Apparently, we have $P(n,m,0)\Rightarrow P(n,m)$. Moreover, we have the following:

\begin{theorem}
Using the above notations, we have
\[
P(n,2,2)\Rightarrow Q(n).
\]
\end{theorem}

{\bf Proof.} Let
\[
X=A_1+A_3,\quad Y=A_2+A_4,
\]
where $A_1, A_2$ are symmetric and $A_3,A_4$ are skew-symmetric matrices. Note that these decomposition are orthogonal.

Using the above notations, we have the decomposition of 
\[
[X,Y]=([A_1,A_4]+[A_3,A_2])+([A_1,A_2]+[A_3,A_4]).
\]
Consequently, we have
\[
||[X,Y]||^2=||[A_1,A_4]+[A_3,A_2]||^2+||[A_1,A_2]+[A_3,A_4]||^2.
\]
Expanding the above expression, we get
\begin{align}\label{long}
\begin{split}
&
||[X,Y]||^2=||[A_1,A_2]||^2+||[A_3,A_4]||^2+||[A_1,A_4]||^2+||[A_2,A_3]||^2\\&
\qquad +2\langle [A_1,A_2],[A_3,A_4]\rangle+2
\langle[A_1,A_4],[A_3,A_2]\rangle.
\end{split}
\end{align}
A straightforward computation gives that
\[
\langle [A_1,A_2],[A_3,A_4]\rangle+
\langle[A_1,A_4],[A_3,A_2]\rangle
=
-\langle[A_1,A_3], [A_2,A_4]\rangle.
\]
Substituting the expression into~\eqref{long} and using the Cauchy inequality, we have
\[
||[X,Y]||^2\leq\sum_{i<j}||[A_i,A_j]||^2.
\]

If $P(2,2,n)$ is true, then we have
\[
||[X,Y]||^2\leq\frac 12(\sum_{i=1}^4||A_i||^2)^2=\frac 12(||X||^2+||Y||^2)^2.
\]

Replacing $X$ by $tX$ and $Y$ by $Y/t$, we have
\[
||[X,Y]||^2\leq\frac 12 (t^2||X||^2+\frac{1}{t^2}||Y||^2)^2.
\]
Minimizing the right hand side of the above with respect to $t$, we get the desired inequality:
\[
||[X,Y]||^2\leq 2||X||^2\cdot||Y||^2.
\]

\qed

In summary, the relations of the conjectures are
as follows:
\[
P(n,m)\Leftarrow P(n,m,0)\Leftarrow P(n,m,m')\Rightarrow P(n,2,2)\Rightarrow Q(n).
\]

\section{Relation to the Pinching theorems}

Let $M$ be an $n$-dimensional compact minimal submanifold in a unit sphere $S^{n+m}$ of dimension $n+m$. Let $||\sigma||^2$ be the square of the length of the second fundamental form. Through the works of Chern-do Carmo-Kobayashi~\cite{chern-d-k}, Yau~\cite{yau-pinching}, Shen~\cite{shen-yibing}, and Wu-Song~\cite{wu-song}, Li-Li~\cite{lianmin} and Chen-Xu~\cite{xusenlin} got the following optimal pinching theorem:

\begin{theorem}
Let $M$ be an $n$-dimensional compact minimal submanifold in $S^{n+m}$, $p\geq 2$. If $||\sigma||^2\leq\frac 23n$ everywhere on $M$, then $M$ is either a totally geodesic submanifold or a Veronese surface in $S^4$.
\end{theorem}

The proof is based on the following type of Bochner formula

\[
\frac 12\Delta||\sigma||^2=\sum_{i,j,k,\alpha}(h^\alpha_{ijk})^2+n||\sigma||^2-\sum_{\alpha,\beta}||[A_\alpha,A_\beta]||^2-\sum_{\alpha,\beta}|\langle A_\alpha,A_\beta\rangle |^2,
\]
where $A_\alpha$ is the matrix $(h_{ij}^\alpha)$; and $h^\alpha_{ij}$ is the second fundamental form with respect to the  orthonormal basis of the tangent
 and the normal bundles; $h_{ijk}^\alpha$ is the covariant derivative of the second fundamental form.

 In ~\cite{lianmin}~\footnote{The proof of ~\cite{xusenlin} is more geometric.}, the following result was proved (cf.~\cite[pp 585, equation (5)]{lianmin}):
 
 \begin{theorem}\label{thm42} Using the same notations as above, we have
  \[
  2\sum_{i<j}||[A_i,A_j]||^2\leq \frac 32\left(\sum_{i=1}^m||A_i||^2\right)^2-\sum_{i=1}^m||A_i||^4.
  \]
  \end{theorem}
  
  \qed
  
  We denote  the above inequality to be $P'(n,m)$.  
  In this section, we prove the following
  
  \begin{theorem}\label{thm}The DDVV conjecture implies Theorem~\ref{thm42}. That is,

  \[
  P(n,m)\Rightarrow P'(n,m).
  \]
  \end{theorem}
  
  Thus inequality $P(n,m)$ is sharper than that in Theorem~\ref{thm42}.
  
  {\bf Proof.} We assume that 
  \[
  ||A_1||\geq\cdots \geq ||A_m||.
  \]
  We prove $P'(n,m)$ by induction: suppose $P'(n,m-1)$ is true. Then we have the following 
  
  \begin{lemma}   If $P'(n,m)$ is true for 
  \[
  ||A_1||^2\leq\sum_{i=2}^m||A_i||^2,
  \]
  then $P'(n,m)$ is true for any $A_1,\cdots,A_m$.
  \end{lemma}
  
  {\bf Proof.} We let $A_1=tA_1'$ and assume that $||A'_1||=1$. Then inequality $P'(n,m)$ can be written as
  \begin{align}\label{b}
  \begin{split}&
  \frac 12 t^4-t^2\left (2\sum_{i=2}^m||[A_1',A_i]||^2-3\sum_{i=2}^m||A_i||^2\right)\\&
  +\frac 32\left(\sum_{i=2}^m||A_i||\right)^2-\sum_{i=2}^m||A_i||^4
  -2\sum_{2\leq i<j}||[A_i,A_j]||^2\geq 0.
  \end{split}
  \end{align}
 By the inductive assumption, the total of the last three terms of the above is nonnegative. Let
 \begin{equation}\label{defa}
  a=2\sum_{i=2}^m||[A_1',A_i]||^2-3\sum_{i=2}^m||A_i||^2.
  \end{equation}
  If $a\leq 0$, then
  then ~\eqref{b} is trivially true. On the other hand, if $a>0$, then the minimum value is obtained at
  \[
  t^2=a.
  \]
Using the fact that $||[A_1',A_i]||^2\leq 2||A_i||^2$ (see \S ~\ref{p1}), we obtain:
  \[
  ||A_1||^2\leq \sum_{i=2}^m||A_i||^2.
  \]
  
  \qed

  {\bf Proof of  Theorem~\ref{thm}.} If 
    \[
  ||A_1||^2\leq\sum_{i=2}^m||A_i||^2,
  \]
  then
  \[
  \left(\sum_{i=1}^m||A_i||^2\right)^2
  \leq   
  \frac 32\left(\sum_{i=1}^m||A_i||^2\right)^2-\sum_{i=1}^m||A_i||^4.
  \]
  Thus
  \[
 P(n,m)\Rightarrow P'(n,m).
   \]

   \qed
  
  \begin{remark}
Since $P(3,m)$, $P(n,3)$ are  true by the  results in~\cite{cl,lu21}, can we get new pinching theorems using these new shaper inequalities?
\end{remark}

We conjecture the following to be true:

\begin{conj} There is a constant $\eps(n)$, depending only on $n$, such that if $||\sigma||^2\leq \frac 23n+\eps(n)$, then $M$ has to be totally geodesic or Veronese surface in $S^4$.
\end{conj}

This conjecture is a more general conjecture of Chern type. See Peng-Terng~\cite{peng-terng}, Cheng-Yang~\cite{cheng-yang} and the references there for details.

 \section{A warming up exercise}\label{p1}
 The following result was proved  in~\cite{ddvv}.
 The proof  is quite easy. However, we go through it  because one can see the difficulties of the DDVV conjecture from the proof.
  
  \begin{theorem}
   $P(n,2)$ is true. That is, if $A,B$ are symmetric matrices. Then
   \[
   ||[A,B]||^2\leq 2||A||^2\cdot||B||^2.
   \]
   \end{theorem}
   
   \begin{remark} As pointed out in~\cite[Theorem 4.1]{bw}, the above inequality is true if one of the matrix is {\it normal}.
   \end{remark}
   
   {\bf Proof.} We let
   \[
   A=Q^T\begin{pmatrix}
   \lambda_1\\&\ddots\\&&\lambda_n\end{pmatrix}
   Q.
   \]
   Let
   \[
   B_1=QBQ^T.
   \]
   Then
   \[
   [A,B]=Q^T[J,B_1]Q.
   \]
   Thus
   \[
   ||[A,B]||^2=||[J,B_1]||^2\leq \sum_{i,j}(\lambda_i-\lambda_j)^2 (b_1)_{ij}^2,
   \]
   where $(b_1)_{ij}$ is the entries of the matrix $B_1$.
   Since 
   \[
   (\lambda_i-\lambda_j)^2\leq 2\sum_i\lambda_i^2,
   \]
   we have
   \[
   ||[A,B]||^2\leq 2\sum_i\lambda_i^2\sum (b_1)_{ij}^2=2||A||^2||B||^2.
   \]
   
   \qed
   
   We make two remarks on the proof.
   \begin{enumerate}
\item The above proof is the ONLY proof we have. It makes use of a non-trivial linear algebraic fact: symmetric matrices are diagnolizable.
   
\item In the proof, we make use of the fact that the Frobienius norm is independent under the change of orthogonal matrices.
   \end{enumerate}
   The conclusion: in order to prove $P(n,m)$, first try to find the invariant group of the inequality. The larger the group, the more reductions (of the inequality)  we can obtain.

\section{Invariance}\label{p0}

Let $A_1,\cdots,A_m$ be  $n\times n$ symmetric matrices.
Let $G=O(n)\times O(m)$. Then $G$ acts on matrices $(A_1,\cdots,A_m)$ in the following natural way: let $(p,q)\in G$, where $p,q$ are $n\times n$ and $m\times m$ orthogonal matrices, respectively. Let $q=\{q_{ij}\}$. Then
\[
(p,0)\cdot (A_1,\cdots,A_m)=(pA_1p^{-1},\cdots, pA_mp^{-1}),
\]
and 
\[
(0,q)\cdot (A_1,\cdots,A_m)=(\sum_{j=1}^m q_{1j}A_j,\cdots,\sum_{j=1}^mq_{mj}A_{j}).
\]

It is easy to verify the following 

\begin{prop}\label{prop1}
Conjecture~\ref{cc2}  is $G$ invariant. That is, in order to prove inequality ~\eqref{cc2-1}
for $(A_1,\cdots,A_m)$, we just need to prove the inequality for any $\gamma\cdot (A_1,\cdots,A_m)$ where $\gamma\in G$.
\end{prop}

\qed

As a consequence of  the above proposition, we have the following interesting

\begin{theorem}\label{thm3}
Let $n\geq 2$ be an integer. If $P(n, \frac 12 n(n-1)+1)$ is true, then $P(n,m)$ is true  for any $m$.
\end{theorem}

{\bf Proof.}  See~\cite{cl}.

\qed

\begin{corollary}\label{cor}
If $P(3,4)$ is true, then $P(3,m)$ is true for $m\geq 2$. \end{corollary}

\qed

\begin{remark}
$G$ is not nearly as big as we expect in the following sense: if we wanted to reduce the problem to the same level as that in~\cite{gmm,guw}, $G$ should have been $O(n)\times O(mn)$.  Usually the smaller the invariant group, the more difficulty the proof of the  inequality will be.
\end{remark}

\section{Sketch of the proofs}

In this section, we sketch of the proof of $P(3,m)$ and $P(n,3)$.

{\bf Proof of $P(3,m)$}. 
Using Corollary~\ref{cor}, we only need  to prove $P(3,4)$. However, the methods of proving $P(3,4)$ and $P(3,3)$ are the same so we only discuss the proof of  $P(3,3)$.
   
We begin with~\cite[\S 4]{cl}:
   
   \begin{theorem}\label{thm71} Let $A,B,C$ be $3\times 3$ symmetric traceless
    matrices. Then
    \[
    \left(||A||^2+||B||^2+||C||^2\right)^2\geq 2||[A,B]||^2+2||[B,C]||^2+2||[C,A]||^2.
    \]
    \end{theorem}
  
  {\bf Sketch of the Proof.}
  The inequality we need to prove contains $15$ independent variables. Our strategy is to reduce the number of independent variables step by step.

   Without loss of generality, we assume that $A$ is diagnolized. Let
  \[
  A=\begin{pmatrix}
  t\eta_1\\& t\eta_2\\&&t\eta_3
  \end{pmatrix}
  \]
  where
  \[
  \eta_1^2+\eta_2^2+\eta_3^2=1.
  \]
  
  Then we have
  \[
  (t^2+||B||^2+||C||^2)^2\geq 2t^2\sum_{i,j}(\eta_i-\eta_j)^2(b_{ij}^2+c_{ij}^2)+2||[B,C]||^2.
  \]
  
  The first reduction: finding the condition such that  the above is true for any $t$. This is doable because the expression is quadratic in $t^2$. 
  
  The second reduction: Maximize the expression
  \[
  \sum_{i,j}(\eta_i-\eta_j)^2(b_{ij}^2+c_{ij}^2)
  \]
  for all $\eta_1^2+\eta_2^2+\eta_3^2=1$. Luckily, for $3\times 3$ matrices, one can get the explicit maximum value.
  
  After the reductions, we get the following inequality: we let 
  \begin{align}
  \begin{split}
 & r_1^2=b_{23}^2+c_{23}^2, r_2^2=b_{13}^2+c_{13}^2,
  r_3^2=b_{12}^2+c_{12}^2;\\
 & |\mu|^2=b_{11}^2+c_{11}^2+b_{22}^2+c_{22}^2+b_{33}^2+c_{33}^2;\\
  &m_0=(r_1^2+r_2^2+r_3^2)^2-3(r_1^2r_2^2+r_2^2r_3^2+r_3^2r_1^2),
  \end{split}
  \end{align}
  Then the inequality is reduced to
  \begin{equation}\label{psq}
  (||B||^2+||C||^2)^2-2||[B,C]||^2\geq 2(\sqrt{m_0}-|\mu|^2)^2
  \end{equation}
  if $2\sqrt{m_0}-|\mu|^2\geq 0$.
  
  Note that using the above two steps, the number of independent variables is reduced to $10$. Of course, the inequality is more nonlinear now.
  
  In \S ~\ref{p1}, we proved that
  \[  (||B||^2+||C||^2)^2-2||[B,C]||^2\geq0.\]
 The inequality ~\eqref{psq} is sharper. The difficulty to prove ~\eqref{psq} is that
  
  \begin{enumerate}
  \item The only way to prove the non-negativeness of the left-hand side of~\eqref{psq} is to diagnolize
  one of the matrices $B$ and $C$, but
  \item $2\sqrt{m_0}-|\mu|^2$ is {\it not} orthogonal group invariant.
  \end{enumerate}
  
  The property comes to help, {\it only} in the $3$-dimensional case, is the following: if

  \[
  \left[\begin{pmatrix}
  0&c&b\\
  c&0&a\\
  b&a&0
  \end{pmatrix},
  \begin{pmatrix}
  0&z&y\\
  z&0&x\\
  y&x&0
  \end{pmatrix}\right]=
  \begin{pmatrix}
  0&r&-q\\
  -r&0&p\\
  q&-p&0
  \end{pmatrix},
  \]
  then
  \[
  \begin{pmatrix}
  p\\q\\r
 \end{pmatrix}=
 \begin{pmatrix}
 a\\b\\c
 \end{pmatrix}\times
 \begin{pmatrix}
 x\\y\\z
 \end{pmatrix}.
 \]
 For the rest of the proof, we first assume that 
 the diagonal part of $B,C$ are zero. Thus we get an inequality of $6$ independent variables. We are able to prove the inequality directly, using the properties of the cross product. Finally, we observe that if the diagonal parts of $B,C$ are not zero, we will get at most a quadratic expression in terms of the diagnol entries of $B,C$. The analysis of the quadratic expression is quite technical and we refer the original paper to the readers.
 
 More recently, in~\cite{lu21}, we remove the assumption in Theorem~\ref{thm71} that the matrices are $3\times 3$. We have
 
 \begin{theorem}\label{thm72} Let $A,B,C$ be $n\times n$ symmetric traceless
    matrices. Then
    \[
    \left(||A||^2+||B||^2+||C||^2\right)^2\geq 2||[A,B]||^2+2||[B,C]||^2+2||[C,A]||^2.
    \]
    \end{theorem}
 
  {\bf Sketch of the Proof.} We need the following two technical lemmas. For the proofs, see~\cite{lu21}.
  
  \begin{lemma} \label{lem1}  Let $x\geq y\geq 0$. Let $(\eta_1,\cdots,\eta_n)$ be  a unit vector. Then if $\{i,j\}\neq\{k,l\}$, we have
\[
(\eta_i-\eta_j)^2 x+(\eta_k-\eta_l)^2 y\leq 2x+y.
\]
\end{lemma}

\qed
\begin{lemma}\label{lem2}
Suppose that $||A||^2+||B||^2+||C||^2=1$ and $||A||\geq||B||\geq||C||$. Let
\[
\lambda={\rm Max}\,(||[A,B]||^2+||[B,C]||^2+||[C,A]||^2),
\]
and let $A,B,C$ be the maximum point. Then we have
\[
2\lambda ||A||^2=||[A,B]||^2+||[A,C]||^2.
\]
\end{lemma}
\qed

  {\bf Continuation of the proof of Theorem~\ref{thm72}.} We assume that
  \[
  ||A||^2+||B||^2+||C||^2=1.
  \]
  Using the above two lemmas, we can get
  \[
  2\lambda||A||^2\leq ||A||^2(2||B||^2+||C||^2)\leq ||A||^2.
  \]
  That is,  $2\lambda\leq 1$, as desired.
  
 \qed
  
  \begin{remark} The same method can be used in the case $m>3$. The details will be in the next version of the paper~\cite{lu21}.
  \end{remark}
  
  \bibliographystyle{abbrv}   
\bibliography{pub,unp,2007}   
\end{document}